\documentclass[11pt]{article}

\usepackage{latexsym}
\usepackage{amsmath}
\usepackage{amsthm}
\usepackage{amssymb}
\usepackage{vmargin}
\usepackage{amscd}
\usepackage{stmaryrd}
\usepackage{euscript}
\usepackage{mathrsfs}
\usepackage{amscd}
\usepackage[all]{xy}
\usepackage{xr}
\DeclareMathAlphabet{\mathpzc}{OT1}{pzc}{m}{it}

\setmargins{32mm}{20mm}{14.6cm}{22cm}{1cm}{1cm}{1cm}{1cm}

\def\ra{\rightarrow}

\def\id{\mathrm{id}}
\def\s|{\,|\,}
\def\ten{\otimes}
\def\vareps{\varepsilon}
\def\eps{\epsilon}

\def\CC{\mathrm{C}}

\def\N{\mathbb{N}}
\def\Z{\mathbb{Z}}
\def\Q{\mathbb{Q}}
\def\R{\mathbb{R}}

\def\Bu{\mathbb{B}}

\def\bB{\mathbb{B}}

\def\bD{\mathbb{D}}

\def\bF{\mathbb{F}}

\def\bH{\mathbb{H}}

\def\bT{\mathbb{T}}

\def\C{\mathcal{C}}
\def\U{\mathcal{U}}

\def\cA{\mathcal{A}}
\def\cB{\mathcal{B}}

\def\cD{\mathcal{D}}
\def\cE{\mathcal{E}}

\def\cH{\mathcal{H}}

\def\cL{\mathcal{L}}

\def\cN{\mathcal{N}}

\def\cS{\mathcal{S}}

\def\cU{\mathcal{U}}

\def\cW{\mathcal{W}}

\def\sA{\mathscr{A}}

\def\sC{\mathscr{C}}

\def\sF{\mathscr{F}}

\def\r{\mathfrak{R}}
\def\Def{\mathfrak{Def}}

\def\H{\mathrm{H}}
\def\z{\mathrm{Z}}

\def\L{\Lambda}

\def\n{\mathfrak{n}}
\def\g{\mathfrak{g}}
\def\ga{\mathpzc{g}}

\def\bu{\ad \Bu}

\def\Hom{\mathrm{Hom}}

\def\Der{\mathrm{Der}}

\def\Iso{\mathrm{Iso}}

\def\Gal{\mathrm{Gal}}

\def\coker{\mathrm{coker\,}}

\def\Ob{\mathrm{Ob}\,}

\def\Set{\mathrm{Set}}

\def\Cat{\mathrm{Cat}}

\def\Mor{\mathrm{Mor}\,}

\def\Grp{\mathrm{Grp}}

\def\Grpd{\mathrm{Grpd}}

\def\Cmpts{\mathrm{Cmpts}\,}
\def\Iso{\mathrm{Iso}}

\def\ad{\mathrm{ad}}

\def\<{\langle}
\def\>{\rangle}
\def\Lim{\varprojlim}
\def\into{\hookrightarrow}
\def\onto{\twoheadrightarrow}
\def\xra{\xrightarrow}

\def\by{\times}

\def\mc{\mathrm{MC}}

\def\ddef{\mathrm{Def}}

\def\mca{\mathrm{MC}_{A}}

\def\ga{\mathrm{G}_{A}}

\def\defa{\mathrm{Def}_{A}}

\def\GL{\mathrm{GL}}

\def\et{\acute{\mathrm{e}}\mathrm{t}}
\def\Et{\acute{\mathrm{E}}\mathrm{t}}

\def\toph{\top_{\mathrm{h}}}
\def\topv{\top_{\mathrm{v}}}
\def\both{\bot_{\mathrm{h}}}
\def\botv{\bot_{\mathrm{v}}}

\def\iff{\Leftrightarrow}

\def\pd{\partial}

\def\half{\frac{1}{2}}

\newtheorem{theorem}{Theorem}[section]

\newtheorem{corollary}[theorem]{Corollary}
\newtheorem{lemma}[theorem]{Lemma}
\newtheorem{theorem*}{Theorem}
\newtheorem{proposition*}[theorem*]{Proposition}
\newtheorem{corollary*}[theorem*]{Corollary}
\newtheorem{lemma*}[theorem*]{Lemma}

\theoremstyle{definition}
\newtheorem{definition}[theorem]{Definition}
\newtheorem{example}[theorem]{Example}

\newtheorem{remark}[theorem]{Remark}
\newtheorem{remarks}[theorem]{Remarks}

\newtheorem{definition*}[theorem*]{Definition}
\newtheorem{example*}[theorem*]{Example}
\newtheorem{examples*}[theorem*]{Examples}
\newtheorem{remark*}[theorem*]{Remark}
\newtheorem{remarks*}[theorem*]{Remarks}

\begin{document}
\title{The structure of the pro-$l$-unipotent fundamental group of a smooth variety\thanks{The author is supported by Trinity College, Cambridge and the Isle of Man Department of Education.}}
\author{J.P.Pridham}
\maketitle

%
\tableofcontents

\section*{Introduction}
\addcontentsline{toc}{section}{Introduction}

In \cite{lgm}, I used the theory of SDCs (developed in \cite{sdc}) to prove results about the structure of the algebraic fundamental group of a smooth variety, similar to those proved in \cite{GM} for the structure of the fundamental group of a K\"ahler manifold. \cite{GM} was motivated by \cite{DGMS}, which shows that the homotopy type of a compact K\"ahler manifold is a formal consequence of its cohomology, the idea being to replace statements about DGAs with those about DGLAs.  

One consequence of \cite{DGMS} is that the de Rham fundamental group of a compact K\"ahler manifold is quadratically presented. In \cite{Morgan}, a weaker result is proved for smooth complex varieties. The primary aim of this paper is to prove an analogous result for smooth proper varieties in finite characteristic. The results in \cite{DGMS} and \cite{Morgan} were inspired by the yoga of weights, motivated by the Weil Conjectures.  With this in mind, it is reasonable to expect the characteristic zero results to have finite characteristic analogues. 

In \cite{Poids} is expressed the philosophy that, for algebraic varieties,  there should be a weight decomposition not only on cohomology,but also on the minimal model (i.e. the homotopy theory), and that only non-zero weights are interesting. In Weil II (\cite{Weil2}), Deligne uses this idea to define the  $\Q_l$-homotopy type of a variety and  prove that it is formal for smooth proper algebraic varieties defined over  finite fields. 

The results of \cite{DGMS} and \cite{Morgan} give explicit information about the structure of the fundamental group. Both papers make use of the Sullivan minimal model, developed in \cite{Sullivan} and shown to determine the de Rham fundamental group $\pi_1(X)\ten \R$ (and, indeed, the entire Whitehead graded Lie algebra). However, since the proof involves the topological methods of \cite{Sullivan}, no such information has been obtained from the formality result in Weil II. By developing a theory of deformations over nilpotent Lie algebras in Section \ref{nilplie}, I obtain similar explicit information on the structure of the algebraic fundamental group. The  Weil Conjectures imply that, for $X$ smooth and proper, the pro-nilpotent fundamental group $\pi_1(X,\bar{x})\ten_{\Z_l}\Q_l$ is quadratically presented. If $X$ is merely smooth, then the pro-nilpotent fundamental group is defined by equations of bracket length at most four. The Hard Lefschetz Theorem implies that, for $X$ smooth and projective, the pro-nilpotent fundamental group cannot be free.

\section{Functors on nilpotent Lie algebras}\label{nilplie}

This section extends the ideas of \cite{Sch} to a slightly different context.

Fix a field $k$, and define $\cN_k$ to be the category of finite-dimensional nilpotent Lie algebras over $k$, and $\widehat{\cN_k}$ to be the category of finitely generated (i.e. $\dim L/[L,L] < \infty$) pro-nilpotent Lie algebras over $k$.

Given $\cL \in \widehat{\cN_k}$, let $\cN_{\cL,k}$ be the category of pairs $(N \in \cN_k, \cL \xra{\phi} N)$, and $\widehat{\cN_{\cL,k}}$ the category of pairs $(N \in \widehat{\cN_k}, \cL \xra{\phi} N)$. We will almost always consider the case $\cL=0$ (note that $\cN_{0,k}=\cN_k$), although a few technical lemmas (Theorems \ref{nMan2} and \ref{nMan3} and Corollary \ref{nkeyhgs}) require the full generality of $\cN_{\cL,k}$. In $\widehat{\cN_{\cL,k}}$, $\cL$ is the initial object, and $0$ the final object.

\begin{definition}\label{nfreelie} Given a finite dimensional vector space $V/k$, denote the free pro-nilpotent Lie algebra on generators $V$ by $L(V)$. Explicitly, let $A(V)$ be the free complete associative algebra on generators $V$,
$$
A(V) := \lim_{\substack{\longleftarrow \\n}} (\bigoplus_m V^{\ten m})/(V)^n,
$$
(the universal enveloping algebra of $L(V)$), with coproduct $\Delta: A(V) \to A(V) \hat{\ten} A(V)$ defined on generators by $\Delta(v)=v \ten 1 + 1 \ten v$. Then
$$
L(V):= \{l \in A(V): \Delta(l)=l \ten 1 + 1 \ten l\}
$$
\end{definition}

\begin{definition}
Given Lie algebras $N,M \in \widehat{\cN_k}$, let $N*M$ be the free Lie algebra product
$$
N*M = \{v \in \U(N)\hat{\ten} \U(M): \Delta(v)=v \ten 1 + 1 \ten v\},
$$
where the coproduct $\Delta$ on $\U(N)\hat{\ten} \U(M)$ comes from those on the (complete) universal enveloping algebras $\U(N)$ and $\U(M)$ (for a definition, see e.g \cite{Se}). Note that $*$ is sum in $\widehat{\cN_k}$ --- the analogue in $\widehat{\C_{\L}}$ is $\hat{\ten}$.
\end{definition}

\begin{definition}
Define the free Lie algebra $L_{\cL}(V):=\cL * L(V)$
\end{definition}

We require that all functors on $\cN_{\cL,k}$ satisfy 
\begin{enumerate}
\item[(H0)]$F(0)=\bullet$, the one-point set.
\end{enumerate}

We adapt the following definitions and results from \cite{Sch} (with identical proofs):

\begin{definition} For $p:N \ra M$ in $\cN_{\cL,k}$ surjective, $p$ is a small extension if \mbox{$\ker p=(t)$,} a principal ideal, such that $[N,(t)]=(0)$. To say that $(t)$ is a principal ideal means that the map
\begin{eqnarray*}
\cL(k) &\to& (t)\\
1 &\mapsto& t
\end{eqnarray*}
is surjective. Note that any surjection can be be factorised as a composition of small extensions.
\end{definition}

For $F:\cN_{\cL,k} \ra \Set$, define $\hat{F}:\widehat{\cN_{\cL,k}} \ra \Set$ by 
$$
\hat{F}(L)= \lim_{\substack{\longleftarrow \\ n}} F(L/\Gamma_n(L)),
$$
 where the $\Gamma_n(L)$ are defined inductively by
$$
\Gamma_1(L)=L, \quad \Gamma_{n+1}(L)=[L,\Gamma_n(L)].
$$
Note that \mbox{$\hat{F}(L) \xrightarrow{\sim} \Hom(h_L,F)$,} where 
\begin{eqnarray*}
h_L:\cN_{\cL,k} &\to& \Set;\\ 
N &\mapsto& \Hom(L,N).
\end{eqnarray*}

\begin{definition}
We will say a functor $F:\cN_{\cL,k} \ra \Set$ is pro-representable if it is isomorphic to $h_L$, for some $L \in \widehat{\cN_{\cL,k}}$. By the above remark, this isomorphism is determined by an element $\xi \in \hat{F}(L)$. We say the pro-couple $(L,\xi)$ pro-represents $F$.
\end{definition}

\begin{definition} A natural transformation $\phi:F \ra G$ in $[\cN_{\cL,k},\Set]$ is called:
\begin{enumerate}
\item unramified if $\phi:t_F \ra t_G$ is injective, where $t_F=F(L(\epsilon))$, where 
$$
L(\eps):=(L(\eps), \cL \xra{0} L(\eps)).
$$
\item smooth if for every $M \twoheadrightarrow N$ in $\cN_{\cL,k}$, we have $F(M) \twoheadrightarrow G(M)\by_{G(N)}F(N)$
\item \'etale if it is smooth and unramified.
\end{enumerate}
\end{definition}

\begin{definition}
$F:\cN_{\cL,k} \ra \Set$ is smooth if and only if $F \ra \bullet$ is smooth.
\end{definition}

\begin{theorem}
\begin{enumerate}
\item Let $A \ra B$ be a morphism in $\widehat{\cN_{\cL,k}}$. Then $h_B \ra h_A$ is smooth if and only if $B=A*L(X_1,\ldots,X_n)$.
\item If $F \ra G$ and $G \ra H$ are smooth morphisms of functors, then the composition $F \ra H$ is smooth.
\item If $u:F \ra G$ and $v:G \ra H$ are morphsms of functors such that u is surjective and $vu$ is smooth, then $v$ is smooth.
\end{enumerate}
\end{theorem}

\begin{definition} A pro-couple $(L,\xi)$ is a hull for $F$ if the induced map $h_L \ra F$ is \'etale.
\end{definition}

\begin{theorem} Let $(L, \xi)$, $(L',\xi')$ be hull of $F$. Then there exists an isomorphism $u:L \ra L'$ such that $F(u)(\xi)=\xi'$.
\end{theorem}

\begin{lemma} Suppose $F$ is a functor such that $$F(V\oplus W) \xrightarrow{\sim}F(V)\by F(W)$$ for vector spaces $V$ and $W$ over $k$, where $V \in \cN_{\cL,k}$ is defined by \mbox{$[V,V]=0$}, \mbox{$\cL \xra{0} V$.} Then $F(V)$, and in particular $t_F$, has a canonical vector space structure, and \mbox{$F(V) \cong t_F \ten V$.}
\end{lemma}

\begin{theorem}\label{nSch}
Given $F:\cN_{\cL,k} \ra \Set$, let $N' \ra N$ and $N'' \ra N$ be morphisms in $\cN_{\cL,k}$, and consider the map:
\begin{enumerate}
\item[$(\dagger)$]  $F(N'\by_N N'') \ra F(N')\by_{F(N)}F(N'').$
\end{enumerate}
Then
\begin{enumerate}
\item $F$ has a hull if and only if $F$ has properties (H1), (H2) and (H3) below:
\begin{enumerate}
\item[{\rm (H1)}] $(\dagger)$ is a surjection whenever $N'' \ra N$ is a small extension.
\item[{\rm(H2)}] $(\dagger)$ is a bijection when $N=0,\quad N''=L(\epsilon)$.
\item[{\rm(H3)}] $\dim_k(t_F) < \infty$.
\end{enumerate}
\item $F$ is pro-representable if and only if $F$ has the additional property (H4):
\begin{enumerate}
\item[{\rm(H4)}]
 $F(N'\by_N N'') \xrightarrow{\sim} F(N')\by_{F(N)}F(N''). $
\end{enumerate}
for any small extension $N' \ra N$.
\end{enumerate}
\begin{proof}
\cite{Sch}, Theorem 2.11.
\end{proof}
\end{theorem}

\begin{definition}
$F:\cN_{\cL,k} \ra \Set$ is homogeneous if 
$$
\eta : F(N'\by_N N'')\ra F(N')\by_{F(N)}F(N'')
$$
is an isomorphism for every $N' \twoheadrightarrow N$.

Note that a homogeneous functor satisfies conditions (H1), (H2) and (H4).
\end{definition}

\begin{definition}
$F:\cN_{\cL,k} \ra \Set$ is a deformation functor if:
\begin{enumerate}
\item $\eta$ is surjective whenever $N' \twoheadrightarrow N$.
\item $\eta$ is an isomorphism whenever $N=0$.
\end{enumerate}
Note that a deformation functor satisfies conditions (H1) and (H2).
\end{definition}

The following results are adapted from\cite{Man}:

\begin{theorem}\label{nSSC} (Standard Smoothness Criterion) Given $\phi :F \ra G$, with \mbox{$(V,v_e) \xrightarrow{\phi'}(W,w_e)$} a compatible morphism of obstruction theories, if $(V,v_e)$ is complete, $V \xrightarrow{\phi'} W$ injective, and $t_F \ra t_G$ surjective, then $\phi$ is smooth.
\begin{proof} \cite{Man}, Proposition 2.17.
\end{proof}
\end{theorem}

For functors $F:\cN_{\cL,k} \ra \Set$ and $G:\cN_{\cL,k} \ra \Grp$, we say that $G$ acts on $F$ if we have a functorial group action $G(N) \by F(N) \xra{*} F(N)$, for each $N$ in $\cN_k$.  The quotient functor $F/G$ is defined by $(F/G)(N)=F(N)/G(N)$.

\begin{theorem}\label{nMan1} If $F:\cN_{\cL,k} \ra \Set$, a deformation functor, and $G:\cN_{\cL,k} \ra \Grp$ a smooth deformation functor, with $G$ acting on $F$, then $D:=F/G$ is a
deformation functor, and if $\nu :t_G \ra t_F$ denotes $h \mapsto h*0$, then $t_D=\coker\nu$, and the universal obstruction theories of $D$ and $F$ are isomorphic.
\begin{proof} \cite{Man}, Lemma 2.20.
\end{proof}
\end{theorem}

\begin{theorem}\label{nMan2} For $F:\cN_{\cL,k} \ra \Set$  homogeneous, and $G:\cN_{\cL,k} \ra \Grp$ a deformation functor, given $a,b \in F(L)$, define $\Iso(a,b): \cN_{L,k} \ra \Set$  by 
$$ 
\Iso(a,b)(L \xrightarrow{f}N)=\{g \in G(N) | g*f(a)=f(b)\}.
$$ 
Then $\Iso(a,b)$ is a deformation functor, with tangent space $\ker\nu$ and, if $G$ is moreover smooth, complete obstruction space $\coker\nu=t_D$.
\begin{proof} \cite{Man}, Proposition 2.21.
\end{proof}
\end{theorem}

\begin{theorem}\label{nMan3} If $G,G'$ smooth deformation functors, acting on $F,F'$ respectively, with $F,F'$ homogeneous, $\ker \nu \ra \ker \nu'$ surjective, and $\coker\nu\ra\coker\nu'$ injective, then $F/G \ra F'/G'$ is injective.
\begin{proof} \cite{Man}, Corollary 2.22.
\end{proof}
\end{theorem}

This final result does not appear in \cite{Man}, but proves extremely useful:

\begin{corollary}\label{nkeyhgs}
If $F:\cN_{\cL,k} \ra \Set$ and $G:\cN_{\cL,k} \ra \Grp$ are deformation functors, with $G$ acting on $F$, let $D:=F/G$, then:
\begin{enumerate}
\item If $G$ is smooth, then $\eta_D$ is surjective for every $M \onto N$ (i.e. $D$ is a deformation functor).
\item If $F$ is homogeneous and $\ker \nu=0$, then $\eta_D$ is injective for every $M \onto N$.
\end{enumerate}
Thus, in particular, $F/G$ will be homogeneous if $F$ is homogeneous, $G$ is a smooth deformation functor and $\ker \nu=0$.\begin{proof}
As for \cite{sdc} Theorem \ref{sdc-keyhgs}.
\end{proof}
\end{corollary}

To summarise the results concerning the pro-representability of the quotient \mbox{$D=F/G$,} we have:
\begin{enumerate}
\item If $F$ is a deformation functor and $G$ a smooth deformation functor, with $\coker \nu$ finite dimensional, then $D$ has a hull.
\item If $F$ is homogeneous and $G$ a smooth deformation functor, with $\ker \nu =0$ and $\coker \nu$ finite dimensional, then $D$ is pro-representable.
\end{enumerate}

\begin{example}
There is one main class of examples of deformation functors on $\cN_k$. Let $k$ be of characteristic zero. Then it is well-known (see e.g. \cite{Am}) that there is an equivalence of categories
$$
\cN_k \xra{\exp} \cU\mathit{ni}_k,
$$
where $\cU\mathit{ni}_k$ is the category of unipotent $k$-groups. Fix a finitely presented group $G$, and consider the functor
\begin{eqnarray*}
h_G:\cU\mathit{ni}_k &\to& \Set\\
U &\mapsto& \Hom(G,U).
\end{eqnarray*}
This is equivalent to the functor
\begin{eqnarray*}
F: \cN_k &\to& \Set\\
N &\mapsto& \Hom(G, \exp(N)).
\end{eqnarray*}
In fact, $F$ is pro-representable by a pro-nilpotent Lie algebra $L$, and by definition 
$$
G\ten_{\Z} k := \exp(L).
$$
\end{example}

\subsection{Differential graded algebras}

\begin{definition} A DGA over a field $k$  of characteristic $0$ is a  graded vector space $A=\bigoplus_{i \in \Z} A^i$ over $k$, equipped with operators $\cdot :A \by A \to A$ bilinear and $d:A \ra A$ linear,  satisfying:

\begin{enumerate}
\item $A^i\cdot L^j \subset A^{i+j}.$

\item $a \cdot b= (-1)^{\bar{a}\bar{b}} b\cdot a$.

\item $(a\cdot b)\cdot c = \cdot (b \cdot c)$.

\item $d(A^i) \subset A^{i+1}$.

\item $d \circ d =0.$

\item $d(a\cdot b) = da\cdot b +(-1)^{\bar{a}} a\cdot db$
\end{enumerate}

Here $\bar{a}$ denotes the degree of $a$, mod $ 2$, for $a$ homogeneous. Note that our DGAs need not have a unit.

\end{definition}

Fix a  DGA $A$.

\begin{definition} The Maurer-Cartan functor 
$\mc_A:\cN_k \ra \Set$ is defined by 
$$ \mc_A(N)=\{x \in A^1\otimes N |dx+\half[x,x]=0\}.$$
\end{definition}

Observe that for $\omega \in A^1\ten N$, 
$$
d\omega +\half[\omega,\omega]=0 \Rightarrow (d+\ad_{\omega})\circ(d+\ad_{\omega})=0,
$$ 
so $(A\ten N, [,], d+\ad_{\omega})$ is a DGLA.

\begin{definition}
Define the gauge functor $\ga:\cN_k \ra \Grp$ by 
$$
\ga(N)=\exp(A^0 \ten N).
$$ 
\end{definition}

We may now define the  DGLA $(A\ten N)_d$ as in \cite{Man}:
$$
(A\ten N)_d^i=\left\{\begin{matrix} (A\ten N)^1\oplus k d	& i=1\\
                           (A\ten N)^i 			& i \ne 1, \end{matrix}\right.
$$
 with 
$$
d_d(d)=0,\quad [d,d]=0,\quad [d,a]_d=da,\quad \forall a \in (A\ten N).
$$ 

\begin{lemma} 
 $exp(A^0 \ten N)$ commutes with $[,]$ when acting  on $(A\ten N)_d$ via the adjoint action.
\end{lemma}

\begin{corollary}
Since $\exp(A^0 \ten N)$ preserves $(A^1\ten N)+d \subset (A\ten N)_d$ under the adjoint action, and 
$$
x \in \mc_A(N) \iff [x+d,x+d]=0,
$$
 the adjoint action of $\exp(A^0 \ten N)$ on $A^1\ten N+d$ induces an action of $\ga(N)$ on $\mc_A(N)$, which we will call the gauge action.
\end{corollary}

\begin{definition} $\defa=\mca/\ga$, the quotient being given by the gauge action \mbox{$\alpha(x)= \ad_{\alpha}(x+d)-d$.}  Observe that $\ga$ is homogeneous. 

 For $a \in \mca(N)$, define $K_a:\cN_{N,k} \to \Grp$ by 
$$
K_a(M)=\exp((d+\ad_a)A^{-1}\ten M).
$$
 Note that this makes sense, since $(d+\ad_a)^2=0$, so 
$$
(d+\ad_a)A^{-1}\ten M \le A^0\ten M
$$
is a Lie subalgebra.
Note that $K_a$ is then a subfunctor of $\Iso(a,a)\le \ga(A)$, so acts on $\Iso(a,b)$ by right multiplication.

Define the deformation groupoid $\mathfrak{Def}_A$ to have objects $\mca$, and morphisms  given by $\Iso(a,b)/K_a$.
\end{definition}

Now, 
$
t_{\ga}=A^0\epsilon ,
$ 
and 
$
t_{\mca}=\z^1(A \epsilon),
$
 with action 
\begin{eqnarray*}
t_{\ga}\by t_{\mca} &\ra& t_{\mca};\\ 
(b,x) &\mapsto& x+db, \text{ so }
\end{eqnarray*}
$$
 t_{\defa}=\H^1(A). 
$$  

\begin{lemma} $\H^2(A)$ is a complete obstruction space for
$\mca$.
\begin{proof} Given a small extension 
$$
e:0 \ra K \ra N \ra M\ra 0,
$$ 
and $x \in \mca(M)$, lift $x$ to $\tilde{x} \in A^1\ten N$, and let 
$$
h=d\tilde{x} +\half [\tilde{x},\tilde{x}] \in A^2\ten N.
$$
 In fact, $h \in A^2\ten K$, as $dx+\half[x,x]=0$. 

Now,
$$
dh=d^2\tilde{x}+[d\tilde{x},\tilde{x}]= [h-\half [\tilde{x},\tilde{x}],\tilde{x}]=[h,\tilde{x}]=0,
$$ 
since $[[\tilde{x},\tilde{x}],\tilde{x}]=0$ and $[K,N]=0$. Let 
$$
v_e(x)=[h]\in \H^2(A\ten K)=\H^2(A)\ten K.
$$
This is well-defined: if $y=\tilde{x}+z$, for $z \in A^1\ten K$, then
$$
dy +\half [y,y]=d\tilde{x}+dz+\half[\tilde{x},\tilde{x}]+\half[z,z]+[\tilde{x},z]=h+dz,
$$ 
as $[K,N]=0$.

This construction is clearly functorial, so it follows that $(\H^2(A),v_e)$ is a complete obstruction theory
for $\mca$.
\end{proof}
\end{lemma}

Now Theorem \ref{nMan1} implies the following:

\begin{theorem} $\defa$ is a deformation functor, $t_{\defa} \cong \H^1(A)$, and $\H^2(A)$ is a complete obstruction theory for $\defa$.
\end{theorem}

The other theorems can be used to prove:

\begin{theorem}\label{nqis}
If $\phi :A \ra B$ is a morphism of DGAs, and 
$$
\H^i(\phi):\H^i(A) \ra \H^i(B)
$$ 
are the induced maps on cohomology, then:
\begin{enumerate}
\item If $\H^1(\phi)$ is bijective, and $\H^2(\phi)$ injective, then $\defa \ra \ddef_B$ is \'etale.
\item If also $\H^0(\phi)$ is surjective, then  $\defa \ra \ddef_B$ is an isomorphism.
\item Provided condition 1 holds, $\mathfrak{Def}_A \to \mathfrak{Def}_B$ is an equivalence of functors of  groupoids if and only if $\H^0(\phi)$ is an isomorphism.
\end{enumerate}
\begin{proof}
\cite{Man}, Theorem 3.1, mutatis mutandis.
\end{proof}
\end{theorem}

\begin{theorem}
If $\H^0(A)=0$, then $\defa$ is homogeneous.
\begin{proof}
This is essentially Theorem \ref{nkeyhgs}, with some allowance made for $K_a$.
\end{proof}
\end{theorem}

Thus, in particular, a quasi-isomorphism of DGAs gives an isomorphism of deformation functors and of deformation groupoids.

\begin{example} Fix a differentiable manifold $X$. The deformation functor
\begin{eqnarray*}
F:\cN_{\R} &\to& \Set\\
\g &\mapsto& \Hom(\pi_1(X,x),\exp(\g))/\exp(\g),
\end{eqnarray*}
where $\exp(\g)$ acts by conjugation, is governed by the DGA
$$
A^{\bullet}:=\Gamma(X,\sA_{\R}^{\bullet}),
$$
the correspondence being as follows. Given $\omega \in \mca$, let
$$
\Bu_{\omega}:=D^{-1}(\omega),
$$
where
\begin{eqnarray*}
D: \exp(\g \ten \sA^0) &\to& \g \ten \sA^1\\
\alpha &\mapsto& d\alpha\cdot\alpha^{-1}.
\end{eqnarray*}
Then $\Bu_{\omega}$ is a principal $\exp(\g)$-sheaf on $X$. There is a monodromy action of $\pi_1(X,x)$ on the stalk $\Bu_{x}$, and this gives us a representation $\pi_1(X,x) \to G$.

That this gives rise to an isomorphism of functors follows from the flabbiness of $\sA^{\bullet}$ (compare tangent and obstruction spaces). Of course, $\defa$ is not pro-representable, but its hull is $\pi_1(X,x)\ten\R$. 
\end{example}

\section{SDCs over $\cN_k$}

The formal definitions, properties and constructions of SDCs over $\cN_k$ hold in exactly the same manner as those defined in \cite{sdc} over $\C_{\L}$. A summary follows.

\begin{definition}\label{nsdcdef} A simplicial deformation complex $S^{\bullet}$ consists of smooth homogeneous functors $S^n:\cN_k \to \Set$ for each $n \ge 0$, together with maps 
$$
\begin{matrix}
\pd^i:S^n \to S^{n+1} & 1\le i \le n\\
\sigma^i:S^{n}\to S^{n-1} &0 \le i <n
\end{matrix},
$$
an associative product $*:S^m \by S^n \to S^{m+n}$, with identity $1: \bullet \to S^0$, where $\bullet$ is the constant functor $\bullet(\g)=\bullet$ on $\cN_k$, such that:
\begin{enumerate}
\item $\pd^j\pd^i=\pd^i\pd^{j-1}\quad i<j$.
\item $\sigma^j\sigma^i=\sigma^i\sigma^{j+1} \quad i \le j$.
\item 
$
\sigma^j\pd^i=\left\{\begin{matrix}
			\pd^i\sigma^{j-1} & i<j \\
			\id		& i=j,\,i=j+1 \\
			\pd^{i-1}\sigma^j & i >j+1
			\end{matrix} \right. .
$
\item $\pd^i(s)*t=\pd^i(s*t)$.
\item $s*\pd^i(t)=\pd^{i+m}(s*t)$, for $s \in S^m$.
\item $\sigma^i(s)*t=\sigma^i(s*t)$.
\item $s*\sigma^i(t)=\sigma^{i+m}(s*t)$, for $s \in S^m$.
\end{enumerate}
\end{definition}

\begin{remark}
If we set $\omega_0$ to be the unique element of $S^1(0)$, then, since $0$ is the initial object in $\cN_k$, we may set
$\pd^0(s)=\omega_0*s$, and $\pd^{n+1}(s)=s*\omega_0$. $S^{\bullet}$ then becomes a cosimplicial complex.
\end{remark}

\begin{definition} Let $t_S^{\bullet}$ be the tangent space of $S^{\bullet}$, i.e. $t_S^n=S^n(L(\eps))$.
\end{definition}

\begin{definition}
Define the Maurer-Cartan functor $\mc_S$ by 
$$
\mc_S(\g)=\{\omega \in S^1(\g) \colon \omega*\omega=\pd^1(\omega)\}.
$$
\end{definition}

\begin{lemma}
$S^0$ is a group under multiplication.
\end{lemma}

Now, if $\omega \in \mc_S(\g)$ and $g \in S^0(\g)$, then $g*\omega*g^{-1}\in \mc_S(\g)$. We may  therefore make the following definition:
\begin{definition}
$$
\ddef_S=\mc_S/S^0,
$$
the quotient being with respect to the adjoint action. The deformation groupoid
$$
\Def_S
$$
has objects $\mc_S$, and morphisms given by $S^0$.
\end{definition}

\begin{lemma}\label{nfaithful}
The action $S^0 \by S^n \to S^n$ is faithful (i.e. $s*h=t$ for some $t$ only if $s=1$).
\end{lemma}
Ths implies that, for all $\omega \in \mc_S(\g)$, $\sigma^0(\omega)=1$.

\begin{definition}
Define the cohomology groups of $S$ to be
$$
\H^i(S):=\H^i(t_S^{\bullet}),
$$
the cohomology groups of the cosimplicial complex $t_S^{\bullet}$.
\end{definition}

\begin{lemma}
The tangent space of $\mc_S$ is $\z^1(t_S^{\bullet})$, with the action of $S^0$ giving $\nu(s)=\pd^1(s)-\pd^0(s)$.
\end{lemma}

\begin{lemma}\label{nobstrsdc}
$\H^2(S)$ is a complete obstruction space for $\mc_S$.
\end{lemma}

Theorems \ref{nSSC} to \ref{nkeyhgs} now imply:

\begin{theorem}
$\ddef_S$ is a deformation functor, with tangent space $\H^1(S)$ and complete obstruction space $\H^2(S)$. Moreover, if $\H^0(S)=0$, then $\ddef_S$ is homogeneous.
\begin{proof}
Theorem \ref{nMan1} and Corollary \ref{nkeyhgs}.
\end{proof}
\end{theorem}

\begin{theorem}\label{nsdcqis}
If $\phi :S \ra T$ is a morphism of SDCs, and 
$$
\H^i(\phi):\H^i(S) \ra \H^i(T)
$$ 
are the induced maps on cohomology, then:
\begin{enumerate}
\item If $\H^1(\phi)$ is bijective, and $\H^2(\phi)$ injective, then $\ddef_S \ra \ddef_T$ is \'etale.
\item If also $\H^0(\phi)$ is surjective, then  $\ddef_S \ra \ddef_T$ is an isomorphism.
\item Provided condition 1 holds, $\mathfrak{Def}_S \to \mathfrak{Def}_T$ is an equivalence of functors of  groupoids if and only if $\H^0(\phi)$ is an isomorphism.
\end{enumerate}
\end{theorem}

Call a morphism $\phi:S \to T$ a quasi-isomorphism if the $\H^i(\phi):\H^i(S) \ra \H^i(T)$ are all isomorphisms.

\begin{definition}\label{ndefmor} Given a morphism $\phi:S \to T$ of SDCs, define the groupoid
$$
\Def_{\phi}
$$
to be the fibre of the morphism
$$
\Def_S \to \Def_T
$$
over the unique point $x_0 \in \mc_T(0)$.

Explicitly, $\Def_{\phi}(\g)$ has objects 
$$
\{(\omega, h) \in \mc_S(\g) \by T^0(\g) \quad\colon\quad h\phi(\omega)h^{-1}=x_0\},
$$
and morphisms 
$$
S^0(\g),\quad\text{ where }\quad  g(\omega,h)=(g\omega g^{-1},h\phi(g)^{-1}).
$$
\end{definition}

\begin{theorem}\label{ncone}
Let $S^n_{\phi}(A)$ be the fibre of $S^n(A) \xra{\phi} T^n(A)$ over $x_0^n$. Then $S^{\bullet}_x$ is an SDC, and the canonical map
$$
\Def_{S_{\phi}} \to \Def_{\phi}
$$
is an equivalence of groupoids.
\end{theorem}

\subsection{Computing SDCs}

Throughout this section, we will consider functors $\cD:\cN_k \to \Cat$. We will \emph{not} require that these functors satisfy (H0) (the condition that $F(0)=\bullet$).

\begin{definition} Given a functor $\cD:\cN_k \to \Cat$, and an object $D \in \Ob\cD(0)$, define
$\Def_{\cD,D}:\cN_k \to \Grpd$ by setting $\Def_{\cD,D}(\g)$ to be the fibre of $\cD(\g) \to \cD(0)$ over $(D,\id)$.
\end{definition}

\begin{definition}
We say a functor $\cE:\cN_k \to \Cat$ has uniformly trivial deformation theory if the functor $\Mor\cE$ is smooth and homogeneous, and the functor $\Cmpts\cE$ is  constant,  i.e for $\g \onto \mathfrak{h}$ in $\cN_k$, $\cE(\g) \to \cE(\mathfrak{h})$ is full and essentially surjective. 
\end{definition}

Assume that we are given functors $\cA,\cB,\cD,\cE:\cN_k \to \Cat$, and a diagram
$$
\xymatrix@=8ex{
\cD \ar@<1ex>[r]^{U}_{\top} \ar@<-1ex>[d]_{V} 
&\ar@<1ex>[l]^{F} \cA  \ar@<-1ex>[d]_{V} 
\\
\ar@<-1ex>[u]_{G}^{\dashv}	\cB \ar@<1ex>[r]^{U}_{\top} 
&\ar@<1ex>[l]^{F} \ar@<-1ex>[u]_{G}^{\dashv} \cE,  
}
$$
where $\cE$ has uniformly trivial deformation theory, the horizontal adjunctions are monadic and the vertical adjunctions comonadic. Let 
\begin{align*}
\toph&=UF&	\both&=FU\\
\botv&=VG&	\topv&=GV,
\end{align*}
with 
$$
\eta:1 \to \toph,\quad \gamma:\botv \to 1,\quad \vareps:\both \to 1 \text{ and }\alpha:1 \to \topv.
$$
Assume moreover that the following identities hold: 
\begin{eqnarray*}
GU&=&UG\\ 
FV&=&VF\\
UV&=&VU
\end{eqnarray*}
\begin{align*}
V\vareps&=\vareps V&	U\alpha&=\alpha U\\
V\eta&=\eta V&			U\gamma&=\gamma U.
\end{align*}

\begin{theorem}\label{nmain}
For $D \in \Ob\cD(0)$, let $E=UV D \in \Ob\cE(0)$, and write $E(\g)$ for the image of $E$ under the map $\Ob\cE(0) \to \Ob \cE(\g)$. Set $S^n(\g)$ to be the fibre
$$
S^n(\g)=\Hom_{\cE(\g)}(\toph^n E(\g), \botv^n E(\g))_{UV(\alpha^n_{D}\circ \vareps^n_D)}
$$
of
$$
\Hom_{\cE(\g)}(\toph^n E(\g), \botv^n E(\g)) \to \Hom_{\cE(0)}(\toph^n E, \botv^n E)
$$
over ${UV(\alpha^n_{D}\circ \vareps^n_D)}$.

We give $S^n$ the structure of an SDC as in \cite{sdc} Theorem \ref{sdc-main}.

Then 
$$
\Def_{\cD,D} \simeq \Def_S.
$$
\begin{proof}
As for \cite{sdc} Theorem \ref{sdc-main}.
\end{proof}
\end{theorem}

\subsection{The pro-unipotent fundamental group}\label{nrep}

\begin{definition}\label{malcev}
Take a  pro-finite group $\Gamma$. Consider the functor on $\cN_{\Q_l}$ given by  
$$
\g \mapsto \Hom(\Gamma, \exp(\g)).
$$
This has tangent space $\H^1(\Gamma, \Q_l)$, and universal obstruction space $\H^2(\Gamma, \Q_l)$. If $\dim \H^1(\Gamma, \Q_l)< \infty$, it is pro-representable, by $\cL\in \widehat{\cN_{\Q_l}}$ say. Define the pro-$\Q_l$-unipotent completion  of $\Gamma$ by
$$
\Gamma \ten_{\hat{\Z}} \Q_l:=\exp(\cL(\Gamma,\Q_l)),
$$
and the $\Q_l$ Malcev Lie algebra by
$$
\cL(\Gamma,\Q_l):=\cL.
$$
\end{definition}

\begin{remark}
It follows that $\cL$ is a hull for the quotient functor 
$$
\g \mapsto \Hom(\Gamma, \exp(\g))/\exp(\g).
$$
\end{remark}

Fix a connected scheme $X$ and a geometric point $\bar{x} \to X$. We wish to study the structure of the group
$$
\pi_1(X,\bar{x})\ten\Q_l.
$$

Consider the functor
$$
\Hom(\pi_1(X,\bar{x}),-):\cN_{\Q_l} \to \Grpd,
$$
which associates to $\g \in \cN_{\Q_l}$ the discrete groupoid
$$
\Hom(\pi_1(X,\bar{x}), \exp(\g)).
$$

Recall the following definition from \cite{lgm} Section \ref{lgm-arbql}:
\begin{definition} Given a pro-$l$ group $K$, define a constructible principal $K$-sheaf to be a principal $K$-sheaf $\bD$, such that 
$$
\bD=\lim_{\substack{\longleftarrow \\{K \to F \text{ finite}}}} F \by^K \bD.
$$
Given an $l$-adic Lie group $G$, a constructible principal $G$-sheaf is a $G$-sheaf $\bB$ for which there exists a constructible principal $K$-sheaf $\bD$, for some $K \le G$ compact, with \mbox{$\bB=G \by^K \bD$} (observe that compact and totally disconnected is equivalent to pro-finite).
\end{definition}

\begin{definition}
Define the functor
$$
\mathfrak{PB}(X):\cN_{\Q_l} \to \Grpd
$$
to associate to each $\g \in \cN_{\Q_l}$, the groupoid of principal constructible $\exp(\g)$-sheaves on $X$.

Define the groupoid 
$$
\mathfrak{PB}_{\bar{x}}(X)(\g)
$$
to be the fibre of
$$
\mathfrak{PB}(X)(\g) \to \mathfrak{PB}(\bar{x})(\g)
$$
over $(\exp(\g),\id)$.
\end{definition}

Alternatively, we may think of $\mathfrak{PB}_{\bar{x}}(X)(\g)$ as
consisting of pairs $(\bB,\theta)$, for 
$$
\bB \in \mathfrak{PB}(X)(\g),
$$ 
and $\theta$ an isomorphism 
$$
\theta:\bB_{\bar{x}} \to \exp(\g).
$$
Thus the crucial difference between the groupoids $\mathfrak{PB}_{\bar{x}}(X)(\g)$ and $\mathfrak{PB}(X)(\g)$ is that only those morphisms which give the identity on $\bB_{\bar{x}}$ are permitted in the former.

\begin{theorem}\label{bunrep}
As in  \cite{sdc} Section \ref{sdc-smoothgroup}, we have a functorial equivalence of groupoids 
$$
\Hom(\pi_1(X,\bar{x}),\exp(\g))/\exp(\g) \to \mathfrak{PB}(X)(\g).
$$

In fact, the same construction gives us an equivalence
$$
\Hom(\pi_1(X,\bar{x}),\exp(\g)) \to \mathfrak{PB}_{\bar{x}}(X)(\g).
$$
\end{theorem}

\begin{remark}
If we have an automorphism  $\sigma$ of the pair $\bar{x} \to X$, it yields automorphisms of both the groupoids above. The construction of the equivalence ensures that it is $\sigma$-equivariant.
\end{remark} 

\begin{definition}
Let 
$$
\mathfrak{FB}(X):\cN_{\Q_l} \to \Cat
$$
associate to $\g$ the category of faithful constructible $\exp(\g)$-sheaves on $X$. Observe that
$\mathfrak{PB}(X)$ is the fibre of $\mathfrak{FB}(X)$ over the trivial $\exp(0)$-sheaf $1$ (the constant sheaf of singleton sets).
\end{definition}

As in \cite{sdc}, we have a functorial comonadic adjunction
$$
\xymatrix{
\mathfrak{FB}(X) \ar@<-1ex>[d]_{u^*}\\
\mathfrak{FB}(X') \ar@<-1ex>[u]_{u_*}^{\dashv},
}
$$
where 
$$
X'=\coprod_{x \in  X} \bar{x},
$$
where for each $x \in X$, a geometric point $\bar{x} \to X$ has been chosen. We have maps $u_x:  \bar{x} \to X$, giving a map $u:X' \to X.$

If $X$ is a variety, we may instead take
$$
X'=\coprod_{x \in  |X|} \bar{x}.
$$

Observe that the second category has uniformly trivial deformation theory, $\exp(\g)$ being smooth, so that deformations are described by the SDC
$$
S^n_X(\g)=\Hom_{\mathfrak{FB}(X')(\g)}(\bT(\g),(u^*u_*)^n \bT(\g)),
$$
where $\bT(\g)=\exp(\g)$, the trivial $\exp(\g)$-sheaf. Note that $(u^*u_*)^n \bT(0)=1$, and $S^n_X(0)$ is one point, so we need not worry about taking the fibre over  $UV(\alpha^n_{D}\circ \vareps^n_D)$ (in the notation of Theorem \ref{nmain}).

Now, write 
$$
\CC^n(X,\sF):= \Gamma(X', (u^{*}u_*)^n u^{*}\sF),
$$
 for sheaves $\sF$ on $X$, and
$$
\sC^n(\sF):=(u_* u^{*})^{n+1}\sF,
$$
so that $\CC^n(X,\sF)=\Gamma(X, \sC^n(\sF))$.

We then  have  isomorphisms
\begin{eqnarray*}
S^n(\g) \cong \Gamma(X',(u^*u_*)^n \bT(\g)) &\cong& \CC^n(X, \exp(\g)),\\
(b \mapsto g\cdot b) &\mapsfrom& g,
\end{eqnarray*}
the latter being a cosimplicial complex of groups, becoming an SDC via the Alexander-Whitney cup product
$$
g*h = (\pd^{m+n}\ldots \pd^{m+2}\pd^{m+1}g)\cdot (\pd^0)^m h,
$$
 for $g \in \CC^m,\,h \in \CC^n$.

\begin{lemma}\label{mittag}
If $\H^i(X,\Z_l)$ and $\H^{i-1}(X,\Z_l)$ are finitely generated, then the cohomology groups of $S^{\bullet}_X$ are isomorphic to 
$$
\H^*(X,\Q_l).
$$ 
\begin{proof}
Observe that the tangent space of $S^{\bullet}$ is
$$
(\lim_{\substack{\longleftarrow \\n}}\CC^{\bullet}(X,\Z/l^n))\ten_{\Z_l}\Q_l,
$$
where
$$
\CC^{\bullet}(X,\Z/l^n))
$$
is a complex computing
$$
\H^*_{\et}(X,\Z/l^n).
$$
Now,
\begin{eqnarray*}
\H^i(S_X)
&=& \H^i( \lim_{\substack{\longleftarrow \\n}} \CC^{\bullet}(X, \Z/l^n)\ten_{\Z_l}\Q_l)\\
&=& \H^i( \lim_{\substack{\longleftarrow \\n}} \CC^{\bullet}(X, \Z/l^n))\ten_{\Z_l}\Q_l\\
&=& \H^i( \lim_{\substack{\longleftarrow \\n}} \CC^{\bullet}(X, \Z_l)\ten \Z/l^n)\ten_{\Z_l}\Q_l,
\end{eqnarray*}
where 
$$
\CC^{\bullet}(X, \Z_l)=\lim_{\substack{\longleftarrow \\n}}\CC^{\bullet}(X, \Z/l^n)),
$$
noting that this is a flat $l$-adic system.

Now, the tower
$$
\cdots \to \CC^{\bullet}(X,\Z_l)\ten \Z/l^{n+1} \to \CC^{\bullet}(X,\Z_l)\ten \Z/l^n \to \cdots
$$
clearly satisfies the Mittag-Leffler condition, so, by \cite{W} Theorem 3.5.8, we have the exact sequence:
$$
\xymatrix@=3ex{
0 \ar[r] & \Lim^1 \H^{i-1}(\CC^{\bullet}(X,\Z_l)_n) \ar[r] \ar@{=}[d] & \H^i(\Lim_n \CC^{\bullet}(X,\Z_l)_n) \ar[r] \ar@{=}[d]&  \Lim_n \H^i(\CC^{\bullet}(X,\Z_l)_n) \ar[r] \ar@{=}[d]&    0\\
0 \ar[r]& \Lim^1 \H^{i-1}(X, \Z/l^n) \ar[r]& \H^i( \Lim_n \CC^{\bullet}(X,\Z_l)_n) \ar[r]& \H^i(X,\Z_l) \ar[r]& 0.
}
$$
From \cite{Mi} Lemma V 1.11, it follows that the $\H^{i-1}(X, \Z/l^n)$ are finitely generated $\Z/l^n$-modules, so satisfy DCC on submodules, so this inverse system satisfies the Mittag-Leffler condition, making the $\Lim^1$ on the left vanish.

Thus
$$
\H^i(S)=\H^i(X,\Z_l)\ten_{\Z_l}\Q_l= \H^i(X,\Q_l).
$$
\end{proof}
\end{lemma}


\begin{lemma}\label{mapschsdc}
Given a map $f:X \to Y$ of schemes, there is a canonical morphism
$$
f^*:S_X^{\bullet} \to S_Y^{\bullet}
$$
of SDCs, inducing the usual map
$$
f^*:\H^i(Y,\Q_l) \to \H^i(X,\Q_l)
$$
on cohomology.
\begin{proof}
The map $f$ induces a map $f':X' \to Y'$, making the following diagram commute:
$$
\begin{CD}
X' @>f'>> Y'\\
@VuVV @VVvV\\
X@>f>> Y.
\end{CD}
$$
We then have
$$
S^n_Y(\g)=\Hom_{Y'}(\g,(v^*v_*)^n\g) \xra{(f')^{\sharp}} \Hom_{X'}((f')^*\g,(f')^*(v^*v_*)^n\g).
$$ 
Now, using the the canonical maps
\begin{eqnarray*}
(f')^*v* &\cong& u^*f^*\\
f^*v_* &\to& u_*(f')^*,
\end{eqnarray*}
we get a map
$$
\Hom_{X'}((f')^*\g,(f')^*(v^*v_*)^n\g)\to \Hom_{X'}((f')^*\g,(u^*u_*)^n(f')^*\g),
$$
but $\g$ is just the constant sheaf, so there is a canonical isomorphism $(f')^*\g \cong \g$. Putting these maps together gives us the morphism
$$
f^*:S_X^{\bullet} \to S_Y^{\bullet}
$$
of SDCs.
\end{proof}
\end{lemma}

\begin{remark}
As in \cite{lgm} Remark \ref{lgm-qlbar}, we may replace $\Q_l$ by any finite extension $E$, or indeed by $\bar{\Q_l}$.
\end{remark}

In fact, we can go further than this in describing $\pi_1(X,\bar{x})\ten\Q_l$:

\begin{theorem}\label{etpointed}
If we consider the morphism
$$
S^{\bullet}_X \xra{u_x^*}  S^{\bullet}_{\bar{x}}
$$
of SDCs, 
define the SDC
$$
S^n_{X,\bar{x}} \subset S^n_X
$$
to be the fibre of $u_x^*$ over the unique $0$-valued point of $S^n_{\bar{x}}$.
Then  
$$
\Def_{S_{X,\bar{x}}}
$$
is equivalent to the discrete groupoid  pro-represented by $\cL(\pi_1(X,\bar{x}),\Q_l)$.

The cohomology of $S^{\bullet}_{X,\bar{x}}$ is
$$
\H^i(X,\bar{x};\Q_l)\cong \left\{\begin{matrix} 0 & i=0 \\ \H^i(X,\Q_l) & i \ne 0, \end{matrix} \right.
$$
provided that the $\H^i(X,\Z_l)$ are finitely generated.

\begin{proof}
Recall, from Theorem \ref{ncone}, that $\Def_{S_{X,\bar{x}}}(\g)$ is equivalent to
$\Def_{u_x^*}(\g)$, the fibre of
$$
\Def_{S_X}(\g)  \xra{u_x^*}  \Def_{S_{\bar{x}}}(\g)
$$
over $(a_0,\id)$, where $a_0 \in S^1_{\bar{x}}(0)$ is the unique $0$-valued object.
From the above discussion, it follows that this groupoid is equivalent to the fibre of
$$
\mathfrak{PB}(X)(\g) \to \mathfrak{PB}(\bar{x})(\g)
$$
over $(\exp(\g),\id)$, which is precisely
$$
\mathfrak{PB}_{\bar{x}}(X)(\g).
$$

Hence, by Theorem \ref{bunrep}, we have a canonical equivalence
$$
\Hom(\pi_1(X,\bar{x}),\g) \to \Def_{S_{X,\bar{x}}}(\g).
$$

The statement about cohomology follows from Theorem \ref{ncone}.
\end{proof}
\end{theorem}

\begin{remark}\label{frobworks}  Observe that, from Lemma \ref{mapschsdc}, an automorphism $\sigma$ of the pair \mbox{$\bar{x} \xra{u_x} X$} will naturally give rise to an automorphism of  
$$
S^{\bullet}_X \xra{u_x^*} S^{\bullet}_{\bar{x}},
$$
and hence of $\Def_{u_x^*}(\g)$. It is clear that the equivalence in the theorem will be $\sigma$-equivariant.
\end{remark}

\begin{remark}\label{tangential} Tangential base points.

If $X=Y-Z$, for $Y$ smooth and $Z$ a closed subscheme, there is a definition of the fundamental group $\pi_1(X,\bar{v})$, where $v$ is a non-zero tangent vector at a point on the boundary of $Z$. This makes use of a morphism of sites:
$$
\Et(X) \to \Et(\bar{v}).
$$
In order to define an SDC governing this fundamental group, we want to replace $X'$ by $X' \sqcup \bar{v}$. Instead, we make use of the morphism
$$
\Et(X) \to \Et(X') \by \Et(\bar{v})= \Et(\bar{v}) \by \prod_{x \in X} \Et(\bar{x}),
$$
which is sufficient. In fact,  working with the sites is implicitly what we were doing before, since infinite coproducts are not defined in the category of schemes.

We use $X'\sqcup \bar{v}$ to define $S^{\bullet}_X$, and the morphism $\bar{v} \to X'\sqcup \bar{v}$ enables the construction of a map $S^{\bullet}_X \to S^{\bullet}_v$, which is what we require.
\end{remark}

\subsection{Equivalence between SDCs and DGAs in characteristic $0$}

The reasoning of \cite{sdc} Section \ref{sdc-sdcdgla} carries over directly to this context. We define the category $\cN_k^{\N_0}$ in the obvious manner.  The only result we need to check is the analogue of \cite{sdc} Lemma \ref{sdc-expolie}:

\begin{lemma}\label{nexpolie}
Given a smooth homogeneous functor $G:\cN_k^{\N_0} \to \Grp$, the tangent space $A^*$ has the structure of a GA (graded-commutative associative algebra without unit). Moreover, there is a canonical isomorphism (functorial in $\g_* \in \cN_k^{\N_0}$):
$$
\exp(\bigoplus_n A^n \ten \g_n) \cong G(\g_*).
$$
\begin{proof}
Since $G$ is smooth and homogeneous, it will be ``pro-represented'' by a ``pro-Artinian'' topological free graded Lie algebra $\cL_*$. The product on $G$ gives rise to a co-associative coproduct
$$
\rho: \cL_* \to \cL_* * \cL_*,
$$
for which $\cL_* \to 0$ is the co-identity (existence of a co-inverse is an automatic consequence of pro-nilpotence).

This induces a map $\rho$ (of Hopf algebras) on the universal enveloping algebras:
$$
\xymatrix{
\cL_* \ar[r]^{\rho} \ar[d] &\cL_* * \cL_* \ar[d]\\
\cU(\cL_*) \ar@{-->}[r]^-{\rho} & \cU(\cL_*)\hat{\ten}\cU(\cL_*).
}
$$
Now,
\begin{eqnarray*}
G(\g_*)&=& \Hom_{\mathrm{Lie}}(\cL^*,\g^*)\\
\bigoplus_n A^n \ten \g_n &=& \Der(\cL^*,\g^*).
\end{eqnarray*}
Hence both embed into $\Hom_{\mathrm{Mod}}(\cL_*, \g_*)$, which is an associative algebra, with product $\alpha\cdot\beta= (\alpha*\beta)\circ \rho$. Here $\alpha*\beta:\cL_* * \cL_* \to \g_*$ arises from the universal property of the free product. On $G(\g_*)$, $\cdot$ is the usual product.

Consider the maps
\begin{eqnarray*}
\log(g) &=&  \sum_{n \ge 1}(-1)^{n-1} \frac{(g-e)^{\cdot n}}{n}\\
\exp(\alpha) &=& \sum_{n \ge 0} \frac{\alpha^{\cdot n}}{n!}
\end{eqnarray*}
on $\Hom_{\mathrm{Mod}}(\cL_*, \g_*)$, where $e:\cL_* \to \g_*$ is the canonical map (i.e. zero).

We may extend Lie algebra homomorphisms (resp. derivations) to homomorphisms (resp. derivations) of the universal enveloping algebras. Since $\cU(\cL_*), \cU(\g_*)$ are associative algebras (with unit), $\exp$ takes derivations to homomorphisms, and $\log$ is its inverse.

Given $\alpha \in \Der(\cL^*,\g^*)=\bigoplus_n A^n \ten \g_n$, let $\alpha'$ be the corresponding derivation 
$$
\alpha':\U(\cL_*) \to \U(\g_*).
$$ 
Then $\exp(\alpha')$ is an algebra homomorphism. But $\alpha=\alpha'|_{\cL_*}$, so $\exp(\alpha)=\exp(\alpha')|_{\cL_*}$, which respects multiplication on $\U(\cL_*)$, so respects the bracket on $\cL_*$. Hence 
$$
\exp(\alpha) \in G(\g_*).
$$ 
Similarly, for $g \in \Hom_{\mathrm{Lie}}(\cL^*,\g^*)=G(\g_*)$, 
$$
\log(g) \in \bigoplus_n A^n \ten \g_n.
$$

It remains only to show that $A^*$ is a graded-commutative algebra. This follows from the functorial (in $\g_*$) isomorphism
$$
\bigoplus_n A^n \ten \g_n = \Der(\cL^*,\g^*),
$$
the right hand object being a Lie algebra, since it is $\log$ of a group.
\end{proof}
\end{lemma}

As in \cite{sdc} Section \ref{sdc-sdcdgla}, we then obtain functors
$$
\xymatrix@1{\mathrm{DGA} \ar@<1ex>[r]^{\cS} & \mathrm{SDC} \ar@<1ex>[l]^{\cD}}
$$
with functorial quasi-isomorphisms $1 \to \cS\cD$, $1 \to \cD\cS$ inducing an equivalence of the localised categories.

\section{Structure of the fundamental group}

\subsection{The topological fundamental group}\label{pi1top}

Fix a connected differentiable manifold $X$, and a point $x \in X$.  Let $\sA^n$ be the sheaf of real-valued $\mathrm{C}^{\infty}$ $n$-forms on $X$.

As in Section \ref{nrep}, the functor
$$
\mathfrak{R}(\pi_1(X,x)):\cN_{\R} \to \Grpd,
$$
which associates to $\g \in \cN_{\R}$ the groupoid with objects
$$
\Hom(\pi_1(X,x), \exp(\g)),
$$
and morphisms $\exp(\g)$ (acting by conjugation), is governed by the SDC
$$
S^n(\g)=\CC^n(X,\exp(\g))= \CC^n(X,\cS^n(\R))(\g).
$$
As in \cite{lgm} Section \ref{lgm-toprep}, we have the following quasi-isomorphisms of SDCs (over $\cN_{\R}$):
$$
\begin{matrix}
\CC^n(X,\cS^n(\R)) \\
\downarrow \\
\CC^n(X,\cS(\sA^{\bullet})^n)\\
\uparrow \\
\cS(\Gamma(X,\sA^{\bullet}))^n,
\end{matrix}
$$
since on cohomology we have:
$$
\begin{matrix}
\H^n(X,\R) \\
\downarrow \\
\bH^n(X,\sA^{\bullet})\\
\uparrow \\
\H^n(\Gamma(X,\sA^{\bullet})),
\end{matrix}
$$
the first quasi-isomorphism arising because $\R \to \sA^{\bullet}$ is a resolution, and
the second  from the flabbiness of $\sA^{\bullet}$.

Hence $\r(\pi_1(X,x))$ is governed by the DGA
$$
A^{\bullet}:=\Gamma(X,\sA^{\bullet}),
$$
whose quasi-isomorphism class, also known as the real homotopy type of $X$, is an invariant which has been studied in great detail.

In \cite{Sullivan}, the minimal model of a connected DGA is defined, from which the de Rham fundamental group $\pi_1(X)\ten \R$ is recovered as the dual Lie algebra. We can, however, recover the de Rham fundamental group as the hull of the functor $\ddef_A$, so do not have to involve minimal models.

We can now adapt the results of \cite{lgm} Section \ref{lgm-toprep} to this context. In fact the proofs are slightly simpler, since we are working with the constant sheaf $\R$, rather than the twisted sheaf $\bu_0$. \cite{lgm} Section \ref{lgm-compactkahler}  then becomes the main result of \cite{DGMS}, \cite{lgm} Section \ref{lgm-smoothcx} becomes the main result of \cite{Morgan}, and the only new result is from adapting \cite{lgm} Section \ref{lgm-arbcx}:

\begin{theorem}
Let $X$ be a complex variety. Then there is a mixed Hodge structure (hence a weight decomposition) on the rational minimal model for $X$. In particular, the Malcev Lie algebra
$$
\cL(\pi_1(X)\ten\R)) \cong \R[[\H^1(X,\R)^{\vee}]]/(f(\H^2(X,\R)^{\vee})),
$$
where $f$ preserves the mixed Hodge structures.
\begin{proof}
Identical to \cite{lgm} Theorem \ref{lgm-arbcxthm}.
\end{proof}
\end{theorem}
This answers a question posed in \cite{Am} as to whether Hodge theory can provide any information about the homotopy type (and hence the structure of the fundamental group) for singular varieties. Since zero weights are permitted in the cohomology of arbitrary varieties, we can, in general, say nothing about the function $f$, vindicating the observation in \cite{Poids} that zero weights are unremarkable. However, if in a particular case $\H^1$ is known to have no component of weight zero, then conclusions can be drawn.

\subsection{The algebraic fundamental group}\label{pi1alg}

Let $k=\bF_q$, take a  connected variety $X_k/k$, and let $X=X_k\ten_k\bar{k}$. Fix a closed point $x$ of $X$, and denote the associated geometric point $x \to X$ by $\bar{x}$. Without loss of generality (increasing $q$ if necessary), we assume that $k(x) \subset \bF_q$.

From Section \ref{nrep}, it follows that the Malcev Lie algebra (see Definition \ref{malcev})
$$
\cL(\pi_1(X,\bar{x}),\Q_l)
$$
pro-represents the groupoid $\Def_{S_{X,\bar{x}}}$, and hence the  functor $\ddef_{S_{X,\bar{x}}}$.

Observe that, since $k(x) \subset \bF_q$, there is a natural action of $\Gal(\bar{\bF_q}/\bF_q)$ on $\pi_1(X,\bar{x})$. Let $F$ be the geometric Frobenius element. It follows from the remarks in Section \ref{nrep} that the action of $F$ on $\cL(\pi_1(X,\bar{x}),\Q_l)$ will correspond to the action of $F^*$ on the  SDC above. 

For $X$ either smooth or proper, \cite{Weil2} shows that all the eigenvalues of Frobenius acting on the cohomology group $\H^i(X,\Q_l)$ are algebraic numbers $\alpha$, and for each $\alpha$, there exists a weight $n$, such that all complex conjugates of $\alpha$ have norm $q^{n/2}$. This provides us with a weight decomposition 
$$
\H^i(X,\Q_l)=\bigoplus_n \cW_n \H^i(X,\Q_l).
$$ 

\begin{theorem}\label{nfrobhull}
There is an isomorphism
$$
\cL(\pi_1(X,\bar{x}),\Q_l) \cong  L(\H^1(X,\Q_l)^{\vee})/(f(\H^2(X,\Q_l)^{\vee})),
$$
where $L(V)$ is the free pro-nilpotent Lie algebra (as in \ref{nfreelie}), and  
$$
f:\H^2(X,\Q_l)^{\vee} \to \Gamma_2L(\H^1(X,\Q_l)^{\vee})
$$ 
preserves the (Frobenius) weight decompositions of \cite{Weil2} 3.3.7. 

Moreover, the quotient map
$$
f:\H^2(X,\Q_l)^{\vee} \to \Gamma_2/\Gamma_3 \cong {\bigwedge}^2(\H^1(X,\Q_l)^{\vee})
$$
is dual to half the  cup product
$$
\H^1(X,\Q_l) \by \H^1(X,\Q_l) \xra{\half\cup} \H^2(X,\Q_l).
$$

\begin{proof}
This will essentially be the same as \cite{lgm} Theorem \ref{lgm-frobhull}. From the remarks in Section \ref{nrep}, it follows that the isomorphism
$$
\Hom(\cL\pi_1(X,\bar{x}),\g) \cong \ddef_{S_{X,\bar{x}}}(\g)
$$
 is Frobenius equivariant. Write 
$$
\cL:=\cL\pi_1(X,\bar{x}).
$$

Since $\cL/\Gamma_n(\cL)$ is a finite dimensional vector space over $\Q_l$ for all $n$, we may use the Jordan decomposition  (over $\bar{\Q_l}$) successively to lift the cotangent space $\cL/\Gamma_2(\cL)$ to $\cL/\Gamma_n(\cL)$, obtaining a space of generators $V \subset \cL$, with the  map
$V \to \cL/\Gamma_2(\cL)$ an isomorphism preserving the weight decompositions. 

Now look at the map 
$$
L(V) \onto \cL;
$$
let its kernel be $J \subset [V,V]$. $(J/[V,J])^{\vee}$ is then a universal obstruction space for $h_{\cL}$. Note that we may lift the Frobenius action on $\cL$ to an action on $L(V)$.

 Since  $\H^2(X, \Q_l)$ is an obstruction space for $\ddef_{\bar{x}}$, there is a unique injective map of obstruction spaces (\cite{Man} Proposition 2.18)
$$
(J/VJ)^{\vee} \into \H^2(X, \Q_l).
$$

Observe that both of these obstruction spaces are Frobenius equivariant. Given  a functor $D$ on which $F$ (Frobenius) acts, we  say an obstruction theory $(W,w_e)$ is Frobenius equivariant if there is a Frobenius action on the space $W$ such that for every small extension
$$
e:\quad 0 \to I \to \mathfrak{h} \to \g \to 0,
$$
and every $d \in D(\g)$, we have $w_e(Fd)=Fw_e(d)$.

Since $h_{\cL} \to \ddef_{\bar{x}}$ commutes with Frobenius, the corresponding map of obstruction spaces does (by uniqueness).

As before, we may use the Jordan decomposition on the successive quotients to obtain a space of generators $W \subset J$, such that $W \to J/[V,J]$ is an isomorphism preserving the  weight decompositions. 

Hence we have  maps
$$
\H^2(X, \Q_l)^{\vee} \onto  J/[V,J] \cong W \into L(V),
$$
and
$$
 V   \cong \cL/[\cL,\cL]  \cong \H^1(X, \Q_l)^{\vee},
$$
preserving the Frobenius decompositions.
Let 
$$
f:\H^2(X, \Q_l)^{\vee} \to L(\H^1(X, \Q_l)^{\vee})
$$
 be the composition.
Then
$$
\cL=L(V)/J \cong L(\H^1(X, \Q_l)^{\vee})/(f(\H^2(X, \Q_l)^{\vee})),
$$
as required.

It only remains to prove the statement about  the cup product:

\begin{lemma} Let $F:\cN_K \to \Set$ be a deformation functor. If $V$ is an obstruction theory for $F$, then there is a canonically defined obstruction map $t_F \by t_F \to V$.
\begin{proof}\emph{(Of Lemma.)}
Since $F$ is a deformation functor, $t_F \by t_F \cong F(L(\eps)\oplus L(\nu))$. Consider the small extension
$$
e:\quad 0 \to \quad K[\eps,\nu] \to L(\eps,\nu)/([\eps,[\eps,\nu]],[\nu,[\nu,\eps]]) \to L(\eps)\oplus L(\nu) \to 0.
$$
Sending $f \in F(L(\eps)\oplus L(\nu))$ to $v_e(f) \in V[\eps,\nu]$ defines the  obstruction map $t_F \by t_F \to V$. Standard arguments show that this is bilinear and anti-symmetric.
\end{proof}
\end{lemma}

Now, the obstruction map
$$
V^{\vee} \ten V^{\vee} \to W^{\vee}
$$
is the dual map to $f: W \to \Gamma_2(V)/\Gamma_3(V)$, while the obstruction map
$$
\H^1(X,\Q_l)\ten \H^1(X,\Q_l) \to \H^2(X,\Q_l)
$$
is the usual cup product. The result now follows from uniqueness of the obstruction map.
\end{proof}
\end{theorem}

\begin{corollary}
If $X$ is smooth and proper, then
$$
\cL(\pi_1(X,\bar{x}),\Q_l)
$$
is quadratically presented. 

In fact,
$$
\cL(\pi_1(X,\bar{x}),\Q_l) \cong  L(\H^1(X,\Q_l)^{\vee})/(\check{\cup}\,\H^2(X,\Q_l)^{\vee})),
$$
where $\check{\cup}$ is dual to the cup product.
\begin{proof}
This follows since, under these hypotheses,  \cite{Weil2} Corollaries 3.3.4--3.3.6 imply that $\H^1(X,\Q_l)$ is pure of weight $1$, and $\H^2(X,\Q_l)$ is pure of weight $2$. 
\end{proof}
\end{corollary}

\begin{example}
This implies that, for $X$ smooth and proper, the pro-$l$ quotient $\pi_1^l(X,\bar{x})$ of $\pi_1(X,\bar{x})$ cannnot be the Heisenberg group 
$$
\cH_3(\Z_l)=\left\{ \begin{pmatrix} 1 & x & y\\0 & 1 & z \\ 0 & 0 & 1 \end{pmatrix} \in \GL_3(\Z_l)\right\},
$$
since this is not of quadratic presentation (in particular, this can be inferred from the non-vanishing of the Massey triple product on $\H^1(\pi_1(X,\bar{x}),\Q_l)$ ---
 see \cite{Am} Ch.3 \S 3 for criteria for a Lie algebra to be quadratically presented).
\end{example}

\begin{corollary}
If $X$ is smooth, then
$$
\cL(\pi_1(X,\bar{x}),\Q_l)
$$
is a quotient of the free Lie algebra $L(\H^1(X,\Q_l)^{\vee})$ by an ideal which is finitely generated by elements of bracket length $2,3,4$.
\begin{proof}
This follows since, under these hypotheses,  \cite{Weil2} Corollaries 3.3.4--3.3.6 imply that $\H^1(X,\Q_l)$ is of weights $1$ and $2$, while $\H^2(X,\Q_l)$ is  of weights $2, 3$ and $4$. Hence the image of $f:\H^2(X,\Q_l)^{\vee} \to L(\H^1(X,\Q_l)^{\vee})$ gives equations:
\begin{eqnarray*}
(\cW_{-1})^2 &\text{weight }-2 \\
(\cW_{-1})^3 +(\cW_{-1})(\cW_{-2}) &\text{weight }-3\\
(\cW_{-1})^4 +(\cW_{-1})^2(\cW_{-2}) +(\cW_{-2})^2 &\text{weight }-4,
\end{eqnarray*}
hence equations of degree at most $4$.
\end{proof}
\end{corollary}

\begin{example}
Thus $\pi_1^l(X,\bar{x})$ cannot be the group
$$
\left\{ \begin{pmatrix} 1 & * & *\\0 & \ddots & * \\ 0 & 0 & 1 \end{pmatrix} \in \GL_5(\Z_l)\right\}.
$$
\end{example}

\begin{remark} If $X$ is singular and proper, weights tell us nothing about the structure of the fundamental group, since zero weights are permitted, so any equations may arise.
\end{remark} 

Finally, a result which does not have an analogue in \cite{lgm}.

\begin{theorem}
Let $X$ be smooth and projective. The cup product
$$
\H^1(\pi_1(X,\bar{x}), \Q_l)\by \H^1(\pi_1(X,\bar{x}), \Q_l) \xra{\cup} \H^2(\pi_1(X,\bar{x}), \Q_l)
$$
is non-degenerate.
\begin{proof}
Recall that $\H^1(\pi_1(X,\bar{x}),\Q_l)$ is a tangent space for our deformation functor, and $\H^2(\pi_1(X,\bar{x}), \Q_l)$ the \emph{universal} obstruction space.
This gives an isomorphism 
$$
\H^1(\pi_1(X,\bar{x}),\Q_l) \xra{\cong} \H^1(X,\Q_l),
$$
and an injection
$$
\H^2(\pi_1(X,\bar{x}),\Q_l) \into \H^2(X,\Q_l),
$$
the maps being compatible with cup products. This yields:
$$
\begin{CD}
\H^1(\pi_1(X,\bar{x}), \Q_l)\ten \H^1(\pi_1(X,\bar{x}), \Q_l) @>\cup>> \H^2(\pi_1(X,\bar{x}), \Q_l)\\
@VV\cong V								@VVV\\
\H^1(X,\Q_l)\ten\H^1(X,\Q_l)					@>\cup>> \H^2(X,\Q_l)\\
@V{\eta^{n-1}\cup}V{\cong}V						@V{\eta^{n-1}\cup}VV	\\				
\H^{2n-1}(X,\Q_l)\ten\H^1(X,\Q_l)					@>\cup>> \H^{2n}(X,\Q_l)\cong \Q_l,
\end{CD}
$$
the isomorphism on the bottom left being the Hard Lefschetz Theorem (\cite{Weil2} Theorem 4.1), $\eta$ being the Chern class of an ample line bundle. The cup product along the bottom is non-degenerate by Poincar\'e duality, hence result.
\end{proof}
\end{theorem}

Since the relations defining the hull arise from the cup product, this indicates that $\pi_1(X,\bar{x})$ must have some defining relations. More precisely:

\begin{corollary}
If $X$ is smooth and proper, then 
$$
\cL(\pi_1(X,\bar{x}))/\Gamma_3(\cL(\pi_1(X,\bar{x}))) \ncong L(V)/\Gamma_3(L(V)),
$$ 
for any free Lie algebra $L(V)$.
\begin{proof}
As for \cite{Am} Proposition 3.25.
\end{proof}
\end{corollary}

\begin{remarks}\label{htpy} Homotopy Theory
\begin{enumerate}

\item
Classically, the DGA $\Gamma(X, \sA^{\bullet})$ (or, rather, its quasi-isomorphism class) is regarded as being the real homotopy type of a differentiable manifold $X$. It would therefore make sense to think of the canonical SDC 
$$
S^n(\g)=\CC^n(X, \exp(\g)),
$$
as capturing the $\Q_l$-homotopy theory of a scheme $X$, since it yields a DGA $\cD(S)$. In fact, $S$ is not merely an SDC, but a cosimplicial complex of group-valued functors. Under the duality of Lemma \ref{nexpolie}, this cosimplicial complex of functors is equivalent to the 
cosimplicial complex
$$
\CC^*(X,\Q_l)
$$
of algebras. It seems reasonable to regard the quasi-isomorphism class of this cosimplicial complex as being the $\Q_l$-homotopy type of $X$. The corresponding differential graded algebra will be quasi-isomorphic to the $\Q_l$-homotopy type defined in \cite{Weil2} \S 5, but the approach used there did not allow any practical consequences to be drawn from the formality of the homotopy type.

If we recall from \cite{Sullivan} that the correct notion of automorphisms of a homotopy type is that of \emph{outer} automorphisms, i.e. automorphisms modulo homotopy, then the reasoning of \cite{Weil2} Corollary 5.3.7 shows that the homotopy type is formal, and that the quasi-isomorphism is Galois equivariant.

\item We may think of the complex $\CC^*(X,\bar{x};\Q_l)$ as being the $\Q_l$-homotopy type of the pair $\bar{x} \to X$. This will be also be formal when $X$ is smooth and proper.
\end{enumerate}
\end{remarks}

\newpage
\bibliographystyle{alpha}
\addcontentsline{toc}{section}{Bibliography}
\bibliography{references.bib}
\end{document}